\documentclass{amsart}

\usepackage{amsmath}
\usepackage{graphicx,color}
\usepackage{amssymb}

\input{epsf.tex}

\usepackage{latexsym}           %%for \Box command
\usepackage{amsfonts}
\usepackage{amsthm}             %%for proof environment
\usepackage{amsmath,amssymb,epsfig}  %% from dissertation

\theoremstyle{definition} %%these items will be in normal text, not italicized
\newtheorem{theorem}{Theorem}[section]
\newtheorem{definition}[theorem]{Definition}
\newtheorem{proposition}[theorem]{Proposition}
\newtheorem{corollary}[theorem]{Corollary}

\begin{document}

\title{The Discrete Fundamental Group of the Order Complex of $B_n$}
\author{H\'el\`ene Barcelo$^*$}

\thanks{$^*$ Supported by an NSA grant, \#H98230-05-1-0256}

\address{Department of Mathematics \& Statistics\\
Arizona State University\\
Tempe, Arizona 85287-1804}
\email{barcelo@asu.edu}

\author{Shelly Smith}
\address{Department of Mathematics \\
Grand Valley State University\\
Grand rapids, Michigan 49401-9403}
\email{smithshe@gvsu.edu}

\thanks{}
\date{\today}

\maketitle

%\today \quad draft  %%this adds the date under the title

%%%%%%%%%%%%%%%%%%%%%%%%%%%%%%%%%%%%%%%%%%%%%%%%%%%%%%%%%%
\thispagestyle{empty}

\begin{abstract}

A few years ago Kramer and Laubenbacher introduced a discrete notion of homotopy for simplicial complexes. 
In this paper, we compute the discrete fundamental group of the order complex of the Boolean lattice.
As it turns out, it is equivalent to computing the discrete homotopy group of the $1$-skeleton of the permutahedron.
To compute this group we introduce combinatorial techniques that we believe will
be helpful in computing discrete fundamental groups of other polytopes.
More precisely, we  use the language of words, over the alphabet of simple transpositions,
to obtain conditions that are necessary and sufficient to characterize the equivalence classes of cycles.
The proof requires only simple combinatorial arguments.
As a corollary, we also obtain a combinatorial proof of the fact that
the first Betti number of the complement of the $3$-equal arrangement is equal to $2^{n-3}(n^2-5n+8)-1.$
This formula was originally obtained by Bj\"orner and Welker in 1995.

\end{abstract}

%%%%%%%%%%%%%%%%%%%%%%%%%%%%%%%%%%%%%%%%%%%%%%%%%%%%%%%%%%
\section{Introduction}
\label{S:Introduction}

In this paper we give an application of a discrete notion of homotopy
theory ($A$-theory hereafter) to subspace arrangements, which yields
an unexpected connection between word problems on the symmetric group
and the computation of the first Betti number of the $k$-equal arrangement
in ${\mathbb R}^n$.
More importantly, in order to compute those Betti numbers we were led to introduce
new combinatorial techniques for evaluating the rank of the discrete fundamental group
of the 1-skeleton of the permutahedron. It is those techniques that we believe will be
useful for computing the discrete fundamental group of the 1-skeleton of other polytopes.

Since the early 1980's, subspace arrangements have been extensively studied
both by topologists and combinatorialists.
One of the reasons for such activities is the fact
that certain complexity theory problems arising in computer science
have a reformulation in terms of subspace arrangements,
thus yielding a beautiful interaction between combinatorics,
topology, geometry, and complexity theory. Bj\"orner gives
a nice account of this interaction in \cite{Bj92} and \cite{BL}.

On the other hand, $A$-theory is
a recently introduced notion of discrete homotopy theory for graphs and
simplicial complexes that was originally developed, in the mid 90's, in the context
of social and technological networks. See \cite{perspectives} for
a survey of this topic, and \cite{BL06}, \cite{HMFB}, \cite{KL}, \cite{malle},
\cite{maurer1}, for applications.
More recently, Dochtermann introduced a different notion
of homotopy of graph maps. His notion is based on the categorical product of graphs
rather than on the cartesian product.  For a nice account of the relation between these
two notions of homotopy see \cite{Do1} and \cite{Do2}.

In 2001, Babson \cite{perspectives} (and independently Bj\"orner \cite{Bjpers})
proved that the discrete fundamental $A$-group, $A_1^{n-3}(\Delta(B_n)),$
of the order complex of the boolean lattice is isomorphic to the fundamental
group, $\pi_1(M_{n,3}),$ of the complement, $M_{n,3},$ of the real
$3$-equal arrangement. Namely,
\begin{eqnarray}
\label{iso}
A_1^{n-3}(\Delta(B_n))\simeq \pi_1(M_{n,3}),
\end{eqnarray}
where $M_{n,3}={\mathbb R}^n - V_{n,3}$ with

$$ V_{n,3}=
\{(x_1,x_2,\ldots ,x_n) \in {\mathbb R}^n |
x_{i_1}=x_{i_2}= x_{i_3},  {\rm  \ for \ some\ }
1\leq i_1 < i_2 <  i_3 \leq n \} .$$

In \cite{BW}, Bj\"orner and Welker showed
that in fact the cohomology groups, $H^i(M_{n,k}),$
of the complement of the $k$-equal arrangements are free. Furthermore, they give
formulae for some of the non-vanishing rank of the cohomology groups.
From these formulae one can deduce that the first Betti number for the
(real) complement of the $3$-equal arrangement is equal to
$$ \text {rank } H^1(M_{n,3})= 2^{n-3}(n^2-5n+8)-1.$$
Their proof is quite intricate and relies on the
Goresky-MacPherson and the Ziegler-\v{Z}ivaljevi\'c
formulae and on some combinatorial methods for computing the homotopy type
of partially ordered sets. A different approach, based on non-pure shellability
of the lattice $\pi_{n,k}$ (the lattice of intersections of the subspaces
associated with the $k$-equal arrangement) can be found in \cite{BjWa}.

Because of the above isomorphism (\ref{iso}), one realizes that computing the
rank of the abelianization of  $A_1^{n-3}(\Delta(B_n))$ will
also yield the first Betti number of $M_{n,3}.$  We show here how to compute
this rank using only combinatorial arguments on $\Gamma(B_n),$ the $1$-skeleton
of the permutahedron, a graph that $A$-theory associates with $\Delta(B_n).$
As it turns out, every discrete homotopy argument can be reformulated
in terms of word problems for the symmetric group $S_n,$ making the arguments
simple and, we hope, shedding new light on the theory of subspace arrangements.

In Section $2,$ we review the fundamental notions of $A$-theory that are needed throughout
this article. In a nutshell, the discrete fundamental group, $A_1^q(\Delta)$
of a simplicial complex $\Delta$ is a certain quotient of its (classical) fundamental
group, $\pi_1(\Delta)$. In order to compute the rank of this group a graph $\Gamma_{\Delta}$
is constructed.  Next, viewing $\Gamma_{\Delta}$ as a $1$-dimensional simplicial complex,
one attaches $2$-cells to all of its $3$ and $4$-cycles, after which
it suffices to compute the fundamental group of this $2$-cell complex.
In practice, one is left to count the number of equivalence classes (under an
appropriate equivalence relation) of cycles of
$\Gamma_{\Delta}$ of length at least 5, which is equal to
$\text{rank\ }A_1^q(\Delta)^{ab}$
when the group is free.

In the case of the order complex of the boolean lattice, $\Delta(B_n),$
the associated graph is $\Gamma (B_n).$
The $n!$ vertices of $\Gamma (B_n)$
can be identified with the permutations of the symmetric group $S_n,$ and
there is an edge between two permutations $\sigma$ and $\tau$ if
there exists a simple transposition, $s_i=(i\ i+1),$ such that $\sigma=\tau s_i.$
Note, we multiply permutations from right to left.
Moreover,
in order to compute the rank of  $A_1^{n-3}(\Delta(B_n))^{ab}$
it will suffice to find the
number of equivalence classes of $6$-cycles in $\Gamma (B_n)$.
Note that in $\Gamma (B_n)$ every {\it primitive} $6$-cycle
(one that is not the concatenation of two $4$-cycles) can be associated with a pair of
consecutive simple transpositions $s_i$ and $s_{i+1}.$
While the permutahedron is a well studied polytope (see for e.g.
\cite{Ziegler}), some of the properties needed
to compute the above mentioned rank come to light when we realize that $\Gamma (B_n)$ can
be obtained by taking the cartesian product of smaller graphs
and then removing some of the edges.  The crux of the argument
relies on the fact that the maximal chains of the direct product of two graded posets
$L_1,\ L_2$ can be expressed as a shuffle of maximal chains from $L_1$ and $L_2.$
The various constructions involved are described in Sections
\ref{S:General Graph} and \ref{S:Graph}. For greater details the interested reader
may wish to consult the second author's Ph.D. thesis \cite{shelly}.

Section \ref{S:Classes} contains the main theorem
of this paper. Namely, (and somewhat informally) two $6$-cycles, $C_1$ and $C_2$ of $\Gamma (B_n)$ associated with
the simple transpositions $s_i$ and $s_{i+1}$ are in the same equivalence
class if and only if  there exists an integer $k\ge 1$ such that
$$ C_2=C_1\tau_1\ldots \tau_k$$
where the $\tau_j$ are transpositions disjoint from  $s_i$ and $s_{i+1}.$
From this result it is relatively easy to compute the total number of equivalence
classes of $6$-cycles, that is
$$\text {rank }  A_1^{n-3}(\Delta(B_n))^{ab} = \text {rank } H^1(M_{n,3}) = 2^{n-3}(n^2-5n +8)-1.$$

In order to prove the main theorem we translate the notion of $G$-homotopic
loops in $\Gamma (B_n)$
to an equivalent notion on a set of words (on the alphabet $\bf{S}$ of all
simple transpositions of $S_n$) that are naturally associated to the loops in
$\Gamma (B_n).$
Two words are equivalent if one can be obtained from
the other by a series of transformations involving only
operations of the form $s_j^2=1$ and $s_ks_j=s_js_k$  for
$|k-j|\geq 2.$ This translation
facilitates the burden of the proof of the main result.
In this paper, all transpositions are simple, and thus
we shall write transposition in lieu of simple transposition.

%%%%%%%%%%%%%%%%%%%%%%%%%%%%%%%%%%%%%%%%%%%%%%%%%%%%%%%%%
\section{Discrete Homotopy Theory for Graphs and Simplicial Complexes}
\label{S:Definitions}

In this section we briefly review some of the basic concepts of discrete homotopy
that will be needed throughout the rest of this paper.
All details can be found in \cite{mainpaper} and \cite{babar}.

\begin{definition}
Let $\Gamma=(V,E)$ and $\Gamma'=(V',E')$ be \textit{simple} graphs, with no loops or parallel edges.
\begin{enumerate}
\item
A \textit{graph map} $f:\Gamma \rightarrow \Gamma'$ is a set map $V \rightarrow V'$
that preserves adjacency, that is, if $vw \in E$, then either  $f(v)$ is
adjacent to $f(w)$ in $\Gamma'$, denoted by $f(v) \sim_{\Gamma'} f(w)$, or $f(v)=f(w)$.

\item
Let $v \in V$ and $v' \in V'$ be distinguished vertices.  A \textit{based graph map}
is a graph map $f:(\Gamma, v) \rightarrow (\Gamma', v')$ such that $f(v)=v'$.
\end{enumerate}
\end{definition}

\noindent We note that if $\Gamma$ is connected, then the image of $f$ is a
connected subgraph of $\Gamma'$.

\begin{definition}
The \textit{cartesian product} $\Gamma \Box \Gamma'$ of two graphs, $\Gamma$ and $\Gamma'$,
is the graph with vertex set $V \times V'$ and an edge between $(v,v')$ and $(w,w')$ if
either
\begin{enumerate}
\item $v=w$ and $v' \sim_{\Gamma'} w'$, or
\item $v \sim_\Gamma w$ and $v'=w'$.
\end{enumerate}
\end{definition}

Let $I_n$ be the path on $n+1$ vertices, with vertices labeled from $0$ to $n$,
and let $I$ be the infinite path with vertices labeled $0, 1,2, \ldots .$
Two based graph maps $f,g: \Gamma \rightarrow \Gamma' $ of simple graphs
are $G$-homotopic (written $f \simeq_G ~ g$)
relative to $v_0'$ and $v_1',$
if there is an integer $n$ and a
graph map $F:\Gamma \Box I_n \rightarrow \Gamma'$ such that

\begin{enumerate}

\item $F(v,0)=f(v)$ and $F(v,n)=g(v)$ for all $v \in V$
\item $F(v_0,j)=v_0'$ and $F(v_1,j)=v_1'$for $0 \leq j \leq n$.
\end{enumerate}

Given a simple graph $\Gamma$
with distinguished vertex $v_0,$ let $A^G_1(\Gamma, v_0)$
be the set of all equivalence classes of based graph maps
$f: I \rightarrow \Gamma$ such that  $f(0)=v_0$ and $f(m)=v_0$
for all $m\geq N_f,$ where  $N_f$ is a positive integer that depends on $f.$
Concatenation of loops is well-defined on this set and it is easy to show
that $A^G_1(\Gamma, v_0)$ is indeed a group.
As in classical topology, if $\Gamma$ is connected, the discrete fundamental group
$A^G_1(\Gamma,v_0)$ is independent of the choice of base vertex.
In this case, we refer to $A^G_1(\Gamma)$ simply as the $A_1$-group of $\Gamma$.

\begin{figure}[ht]
\begin{center}
\includegraphics[height=3in]{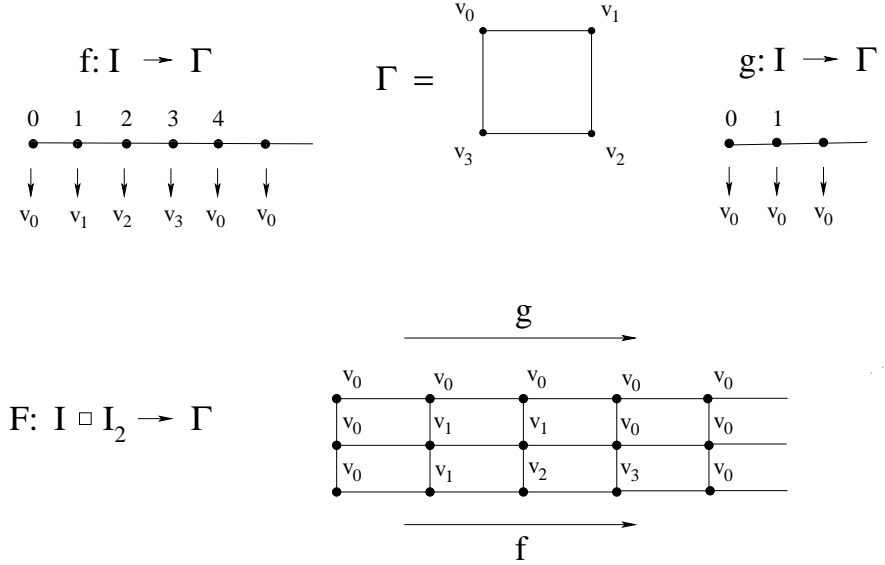}
\caption{A $G$-homotopy from $f$ to $g$.}
\label{F:homotopyFtoG}
\end{center}
\end{figure}

Figure \ref{F:homotopyFtoG} shows an example of two $G$-homotopic graph maps,
$f,g : I \rightarrow \Gamma$ where $\Gamma$ is a cycle of length $4.$
Graph map $f$ corresponds to going around the 4-cycle once, while $g$ is the constant map
$v_0$.  The vertices of the graph (grid) $I \Box I_2$ are labeled with
the image of a $G$-homotopy from $f$ to $g$. The $G$-homotopy $F$ is itself a graph map,
and as such must preserve all adjacencies. Thus for each horizontal or vertical
edge $(v_i,v_j)$
in the grid, we must have either $(v_i,v_j)$ is an edge of $\Gamma$,
or $v_i=v_j$.

Furthermore, it is straightforward to show that any based graph map from $I$ to
the 4-cycle is $G$-homotopic to the constant map $g$, so the $A_1$-group of
the 4-cycle, and similarly of the 3-cycle, is trivial.   However,
a graph map that maps onto a cycle of length $\geq 5$
is not $G$-homotopic to the constant map.
In fact, it can be shown that  if $\Gamma$ is a cycle of length at least $5$, then
$A^G_1(\Gamma)\simeq \mathbb{Z}$ .

In \cite{mainpaper}, Barcelo et al. also show that
$A^G_1(\Gamma,v) \simeq \pi_1(\Gamma,v)/N$,
where $\pi_1(\Gamma,v)$ is the classical fundamental group of $\Gamma$ when
viewed as a 1-dimensional simplicial complex and $N$ is the normal subgroup generated
by 3-and 4-cycles.  Thus, computing the $A^G_1$-group of a graph is equivalent to
attaching 2-cells to the 3- and 4-cycles of the graph and computing the classical
fundamental group of the resulting $2$-cell complex.

\begin{figure}[ht]
\begin{center}
\scalebox{.30}{\includegraphics[height=4in]{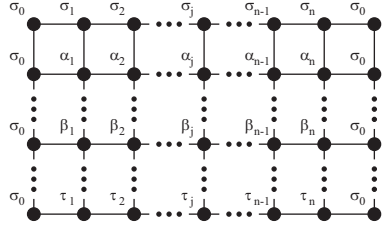}}
\caption{Equivalent $q$-loops.}
\label{F:grid}
\end{center}
\end{figure}

There is an equivalent definition of discrete homotopy for simplicial complexes which
includes a graded version of the discrete fundamental group, related to the
dimension of the intersection of simplices.
Let $\Delta$ be a simplicial complex of dimension $d,$ let $0\leq q\leq d$ be
a fixed integer, and let $\sigma_0\in\Delta$ be a
maximal simplex (with respect to inclusion) of dimension at least $q$.
A $q$-chain in $\Delta$ is a sequence
of simplices (not necessarily distinct),
$$
\sigma,\sigma_1,\sigma_2,\ldots ,\sigma_n,\tau,
$$
such that any two consecutive simplices share a $q$-face.
A $q$-{\it loop in} $\Delta$ {\it based at} $\sigma_0$ is
a $q$-chain beginning and ending at $\sigma_0$.
Two such loops are equivalent if they can be deformed into each other
without breaking any $q$-dimensional connections.
More precisely, the equivalence relation, $\simeq_A,$ on the collection of $q$-loops
in $\Delta$, based at $\sigma_0$, is generated by the following two
conditions.
\begin{enumerate}
\item
The $q$-loop
$$
(\sigma)=(\sigma_0,\ldots ,\sigma_i,\sigma_{i+1},\ldots ,\sigma_n,\sigma_0)
$$
is equivalent to the $q$-loop
$$
(\sigma')=(\sigma_0,\ldots ,\sigma_i,\sigma_i,
\sigma_{i+1},\ldots,\sigma_n,\sigma_0).
$$
That is, we can ``stretch'' loops by repeating a simplex
without changing its equivalence class.
\item
Suppose that $(\sigma)$ and $(\tau)$ have the same length.
They are equivalent if there is a diagram as in Figure \ref{F:grid}.
The diagram is to be interpreted as follows.
A horizontal or vertical edge between two
simplices indicates that they share a $q$-face.
Each horizontal row in the diagram is a $q$-loop based at $\sigma_0$.
Thus, $(\sigma)$ is equivalent to $(\tau)$
if there is a sequence of $q$-loops based at $\sigma_0$ connecting them.
\end{enumerate}

This equivalence relation is called $A$-{\it homotopy}, and
its ensuing set of equivalence classes is denoted by $A_1^q(\Delta,\sigma_0)$.
Concatenation of $q$-loops yields a group structure on
$A_1^q(\Delta,\sigma_0)$. In \cite{mainpaper}, it was shown that in fact
$A_1^q(\Delta,\sigma_0)\cong \pi_1(\Gamma^q_{\Delta},v_0)/N,$ where
$\Gamma^q_{\Delta}$ is the graph whose vertices correspond
to all maximal simplices of $\Delta$ of dimension at least
$q$.  Two vertices $v_{\sigma}$ and $w_{\tau}$ in $\Gamma^q_{\Delta}$ are adjacent
if and only if the corresponding simplices
$\sigma$ and $\tau$ share (at least) a $q$-face, and $v_0$
is the distinguished vertex of $\Gamma^q_{\Delta}$ corresponding
to $\sigma_0$. One sees that there is a close relation between
the $A^G_1$-groups defined for graphs, and the $A^q_1$-groups
defined for simplicial complexes.
Indeed the relation is given by
$$A_1^q(\Delta, \sigma_0) \cong A^G_1(\Gamma^q_{\Delta},v_0),$$
thus $G$-homotopy and $A$-homotopy are equivalent concepts.

%%%%%%%%%%%%%%%%%%%%%%%%%%%%%%%%%%%%%%%%%%%%%%%%%%%%%%%%%
\section{The Product of Lattices}
\label{S:General Graph}

Recall that one of our goals is to use the techniques of $A$-theory to compute the
first Betti number of $M_{n,3}, $ the complement of the $3$-equal arrangement.
As mentioned in the introduction, in \cite{perspectives}
it was shown that
$$A_1^{n-3}(\Delta(B_n))\simeq \pi_1(M_{n,3}),$$
where $\Delta(B_n)$ is the order complex of the boolean lattice,
$B_n-\{\hat 0, \hat 1\}$.
Thus to find $\text {rank } H^1(M_{n,3})$ we will need to count the number of
distinct equivalence classes of cycles, $[C],$ in the graph $\Gamma^{n-3}_{\Delta(B_n)}.$
From here on, for any poset $L$ of rank $r,$ we will only
be interested in its top $A^q_1(\Delta(L))$,
which is $A^{r-3}_1(\Delta(L))$, and thus we shall
write $\Gamma (L)$ in lieu of $\Gamma^q_{\Delta(L)}.$
The vertices of $\Gamma(B_n)$ correspond to the maximal faces of $\Delta(B_n)$,
which are the maximal chains in $B_n-\{\hat 0, \hat 1\}$, which further
correspond to the permutations of $S_n$.  Two vertices in $\Gamma(B_n)$ are adjacent
if the two maximal chains $C_1$ and $C_2$ differ in precisely one element, in which
case we also say that the chains are adjacent, $C_1\sim C_2.$  Equivalently, the
associated permutations differ by multiplication on the right by
a (simple) transposition $s_i$ for some $1 \leq i \leq n-1$.

\begin{figure}[ht]
\begin{center}
\includegraphics[height=1.7in]{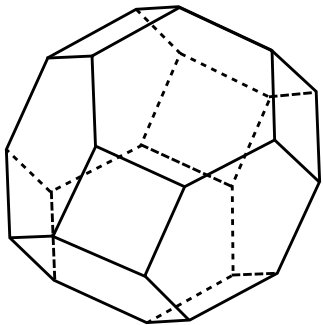}
\caption{The 1-skeleton of the permutahedron $P_3$.}
\label{F:permutahedron}
\end{center}
\end{figure}

We note that $\Gamma(B_n)$ is the 1-skeleton of the permutahedron
$P_{n-1}$ \cite{Ziegler}; Figure \ref{F:permutahedron} represents $\Gamma(B_4)$.
Any path in $\Gamma(B_n)$ corresponds to a product of transpositions,
$s_{i_1}, \ldots s_{i_p}=w$. We view $w$ both as a word in the alphabet $\bf{S},$
as well as an element of the symmetric group $S_n.$
Thus any cycle is a representation of the identity as a product of transpositions.
We can see that if we attach 2-cells to the 4-cycles in $\Gamma(B_4)$,
we are left with eight 6-cycles.  In general, in order to compute the rank of $A^G_1(\Gamma(B_n))^{ab}$, we will need to count the distinct equivalence
classes of cycles, with primitive $6$-cycles as representatives.
Moreover, $A^G_1(\Gamma(B_4))^{ab}$ is a free group (see \cite{perspectives})
on seven generators, not eight, so
unfortunately simply counting 6-cycles will not suffice.

The breakthrough that allows us to understand the $G$-homotopy relation on
$\Gamma(B_n)$ is the simple observation that $B_n$ can be viewed as the direct
product of smaller boolean lattices; in fact, $B_n \simeq B_1^n$.  It is useful
to express this isomorphism as $B_n \simeq B_{n-1} \times B_1$, because the
graph $\Gamma(B_1)$ is a single vertex corresponding to the empty chain.
However, clearly $\Gamma(B_n)\not\simeq \Gamma(B_{n-1}) \Box \Gamma(B_1)$,
since $\Gamma(B_n)$ has $n!$ vertices, compared to $(n-1)!$ vertices in
$\Gamma(B_{n-1}) \Box \Gamma(B_1)$.
In general,
$\Gamma(L_1 \times L_2)\not\simeq \Gamma(L_1) \Box \Gamma(L_2)$
for nontrivial posets $L_1$ and $L_2$, nevertheless, there is a relationship between
the graphs. We now introduce a method to obtain $\Gamma(L_1 \times L_2)$
from $\Gamma(L_1)$ and $\Gamma(L_2)$.

Each maximal chain in $L_1 \times L_2$,
may be viewed as a combination of one
maximal chain from each of $L_1$ and $L_2$; however, those two chains can be
combined in more than one way (for more details on product of posets, see \cite{stanley1}).
Thus, $\Gamma(L_1) \Box \Gamma(L_2)$ is
a subgraph of $\Gamma(L_1 \times L_2)$, but we must also find a way to
reflect the various combinations of chains that are possible.

To remedy this problem, we first introduce a new graph,
\textit{the shuffle graph}, with vertices corresponding to each possible
shuffle of a pair of maximal chains from $L_1$ and $L_2$.  Next, we
construct the cartesian product of the shuffle graph with $\Gamma(L_1) \Box \Gamma(L_2)$.
This solves the problem of too few vertices, but replaces it with a new obstacle
of too many edges.  Finally, we determine which edges are superfluous and remove them so
that the resulting graph is the desired $\Gamma(L_1 \times L_2)$.  In the following
section, we then apply this construction to $B_{n-1} \times B_1$ to create a
representation of the permutahedron, $\Gamma(B_n),$
where we can better understand its structure.

\vskip .2in

\begin{center}
\large\textbf{Step 1 \quad The shuffle graph, $\Gamma^{k,l}_{shuffle}$.}
\end{center}

To construct $\Gamma(L_1 \times L_2)$, we begin by considering the ways in which the edges of chains from $L_1$ and $L_2$ may be combined to create a new chain.
Let $C_1$ and $C_2$ be maximal chains in two graded posets $L_1$ and $L_2$ of
rank $k$ and $l$, respectively.  A shuffle of the edges of $C_1$ and $C_2$ creates a
maximal chain $C$ in $L_1 \times L_2$.  In $C$, label each edge from $C_1$ with the
number of edges from $C_2$ below it in the shuffle.  The ordered, weakly
increasing collection of labels, $\kappa=(a_1,a_2, \ldots ,a_k),$ is the
\textit{$k$-sequence} associated with
that shuffle.  Similarly, label each edge from $C_2$ with the number of edges
from $C_1$ below it in the shuffle and the ordered collection of labels is
an \textit{$l$-sequence}, $\lambda$.

\begin{figure}[ht]
\begin{center}
\includegraphics[height=2.1in]{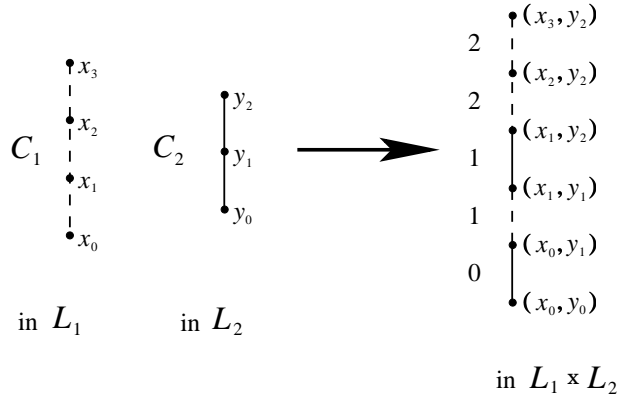}
\caption{The shuffle of $C_1$ and $C_2$ associated with the 3-sequence $(1,2,2)$ and 2-sequence $(0,1)$.}
\label{F:Shuffle}
\end{center}
\end{figure}

We now introduce the shuffle graph, $\Gamma^{k,l}_{shuffle},$  from which we will
construct $\Gamma(L_1 \times L_2).$
The vertices of $\Gamma^{k,l}_{shuffle}$ correspond to
the $\binom{k+l}{k}$ shuffles of chains of $L_1$ and $L_2$
of length $k$ and $l$ respectively.  Label each vertex
with the pair ($\kappa$, $\lambda$) that corresponds to each shuffle.  A
shuffle is uniquely determined by either its $k$-sequence or $l$-sequence,
but both will be useful later in the construction of $\Gamma(L_1 \times L_2)$.

\begin{figure}
\begin{center}
\includegraphics[height=2 in]{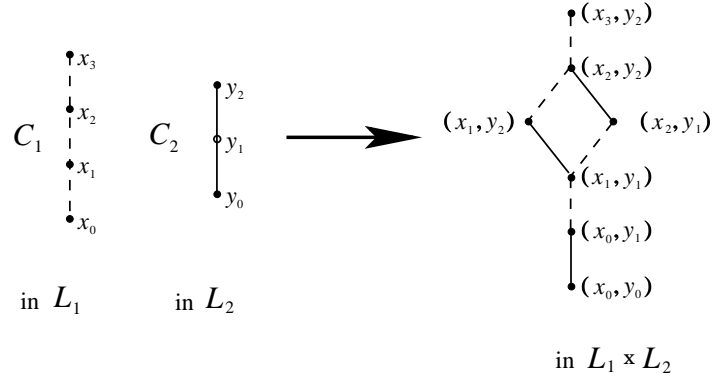}
\caption{Chains $C_1$ and $C_2$ combined using adjacent shuffles with 3-sequences $(1,2,2)$ and $(1,1,2)$.}
\label{F:type1}
\end{center}
\end{figure}

\begin{definition}

Two $k$-sequences, $\kappa=(a_1, \ldots a_k)$ and
$\kappa'=(a'_1, \ldots a'_k)$ are \textit{adjacent,\ }$\kappa \sim \kappa',$ if and only if

\begin{enumerate}
\item $a_i=a'_i\pm 1$ for some $1\leq i \leq k,$ and
\item $a_j=a'_j\ \forall j\neq i.$
\end{enumerate}
Two shuffles of $\Gamma^{k,l}_{shuffle}$ are said to be
adjacent if their $k$-sequences are adjacent, $\kappa\sim\kappa'.$

\end{definition}

\begin{figure}
\begin{center}
\includegraphics[height=2.3in]{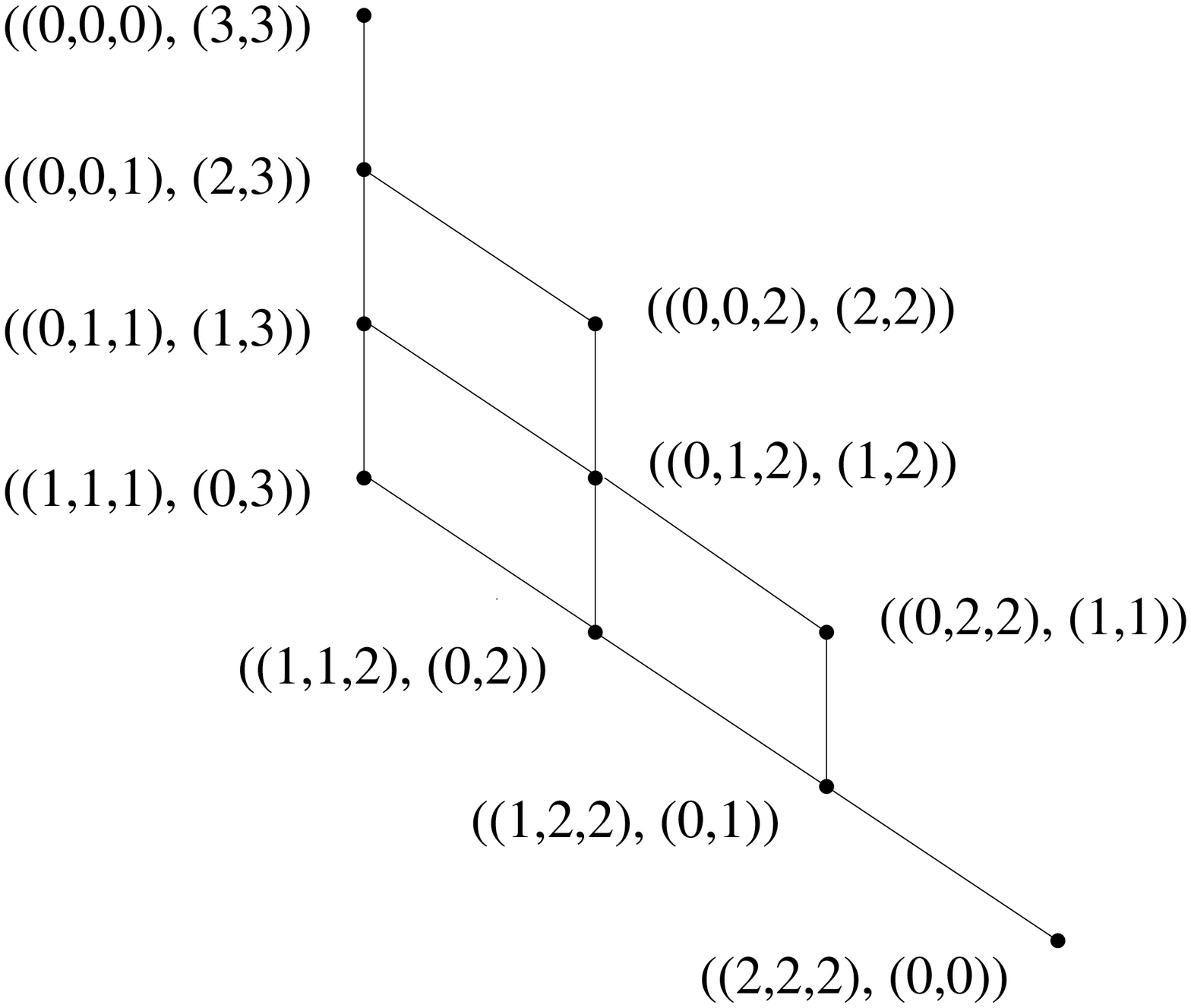}
\caption{$\Gamma_{shuffle}^{3,2}$ labeled with 3-sequences and 2-sequences.}
\label{F:ShuffleGraph}
\end{center}
\end{figure}

We note that if two $k$-sequences are adjacent then
the associated $l$-sequences are also adjacent, so we only
need to refer to one of the sequences when determining if two shuffles
are adjacent.  A pair of chains in $L_1 \times L_2$ resulting from the
use of adjacent shuffles differ by a diamond (formed where the order of the pair
of edges is reversed), as shown in Figure \ref{F:type1}.
Figure
\ref{F:ShuffleGraph} shows $\Gamma^{3,2}_{shuffle}$ with sequences for the
ten possible shuffles of $C_1$ and $C_2$ from Figure \ref{F:type1}.

\vskip .2in

\begin{center}
\large \textbf{Step 2 \quad The intermediate graph $\mathbf{\widetilde{\Gamma}(L_1 \times L_2)}$.}
\end{center}

Let $\widetilde{\Gamma}(L_1 \times L_2)=\Gamma(L_1)\Box \Gamma(L_2) \Box
\Gamma^{k,l}_{shuffle}.$
Label each vertex of
$\widetilde{\Gamma}(L_1 \times L_2)$ with the ordered triple
$(C_1, C_2, (\kappa, \lambda )),$ for all maximal chains $C_i\in L_i,$ $i=1,2$ and
for all possible shuffles $(\kappa, \lambda).$ The
set of vertices of $\widetilde{\Gamma}(L_1 \times L_2)$ corresponds to
all possible shuffles of pairs of maximal chains from $L_1$ and $L_2$,
thus there is a one-to-one correspondence between the vertices of
$\widetilde{\Gamma}(L_1 \times L_2)$ and the maximal chains of $L_1 \times L_2$.
From the definition of a cartesian product of graphs, two vertices in
$\widetilde{\Gamma}(L_1 \times L_2)$,
$(C_1, C_2, (\kappa, \lambda ))$ and
$(C'_1, C'_2, (\kappa', \lambda' )),$ are
adjacent if they satisfy precisely one of the following conditions:

\begin{enumerate}
\item
$C_1=C_1'$, $ C_2=C_2'$, and $\kappa\sim\kappa'$
\item
$C_1=C_1'$, $C_2 \sim C_2'$ in $L_2$, and $\kappa=\kappa'$
\item
$C_1 \sim C_1'$ in $L_1$, $C_2=C_2'$, and $\kappa=\kappa'.$
\end{enumerate}

\begin{figure}[h]
\begin{center}
\includegraphics[height=1.8in]{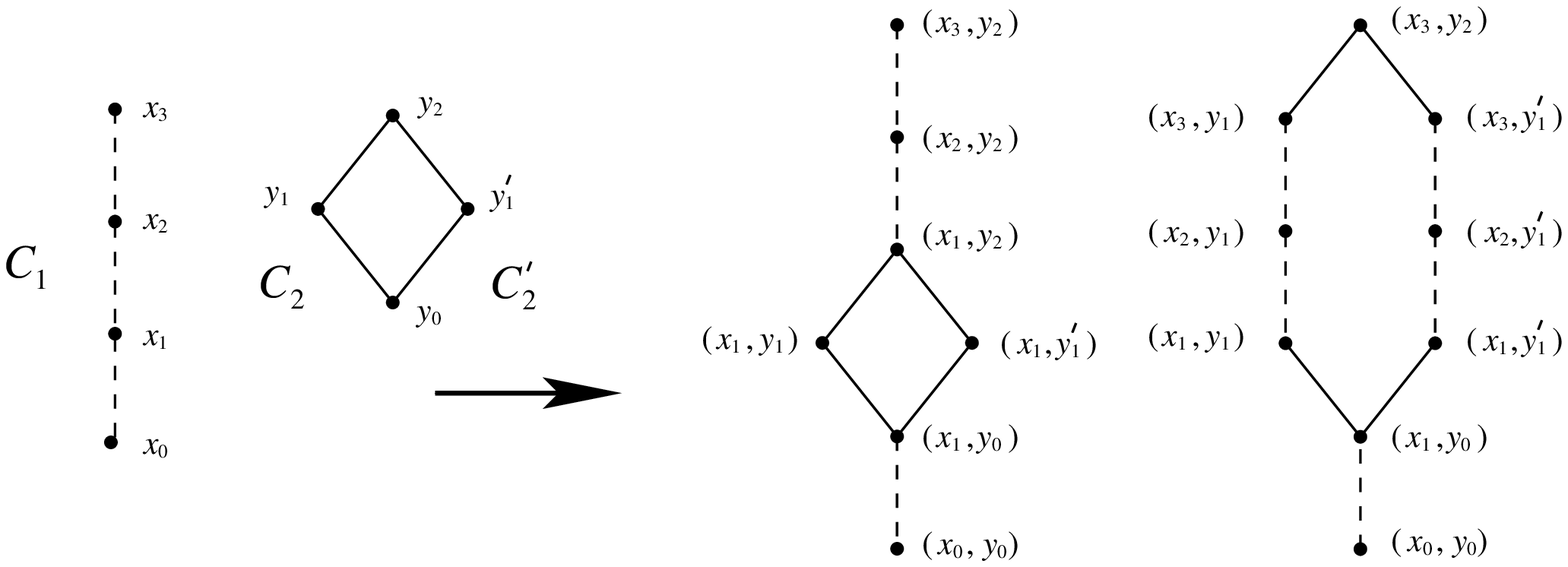}
\caption{Combining $C_1$ with $C_2$ and $C_2'$ using two different shuffles.}
\label{F:type2}
\end{center}
\end{figure}

While $\widetilde {\Gamma}(L_1\times L_2)$ has the right number of vertices,
it has too many edges.
For example, Figure \ref{F:type2} shows the possible result of shuffling
$C_1$ with adjacent chains $C_2$ and $C_2'.$
In one case, the shuffle results in a pair of adjacent chains in
$L_1 \times L_2$.  However, using another shuffle results in chains that differ
at three ranks in $L_1 \times L_2$.  This difficulty is easily remedied in
Step 3 by removing a well-defined set of edges, determined by the rank where
a pair of adjacent chains in one poset differs, along with the $k$- and
$l$-sequences of the shuffle used.

\begin{center}
\large \textbf{Step 3 \quad Removing edges from $\mathbf{\widetilde{\Gamma}(L_1 \times L_2)}$.}
\end{center}

Each edge in $\widetilde{\Gamma}(L_1 \times L_2)$ can be classified as
Type 1, 2, or 3, according to which of the above conditions is
satisfied by $(C_1, C_2, (\kappa, \lambda))$
and $(C'_1, C'_2, (\kappa', \lambda')).$ The
next step in the process of constructing $\Gamma(L_1 \times L_2)$ is to examine
each type of edge in $\widetilde{\Gamma}(L_1 \times L_2),$
to determine which ones are between a pair of adjacent chains
in $L_1 \times L_2$ and which are not.  Once edges corresponding to pairs of
non-adjacent chains have been removed from the graph, the result will be the
desired final graph $\Gamma(L_1 \times L_2)$.

\textbf{Type 1 edges}. $C_1=C_1'$, $ C_2=C_2'$, and $\kappa \sim \kappa'.$
It is easy to see that none of these edges need to be removed.
See Figure \ref{F:type1} for an example of this type.

\textbf{Type 2 edges}. $C_1=C_1'$, $C_2 \sim C_2'$ in $L_2$, and
$\kappa = \kappa'.$
The diagram of $C_2$ and $C_2'$ contains a diamond in $L_2$ at some rank $i$
where the two chains differ.  When we shuffle $C_2$ and $C_2'$ with $C_1$,
this diamond may be \textit{stretched} by the insertion of edges from $C_1$,
depending on which shuffle is used. If $i\notin \kappa$ then the resulting chains
are adjacent; but if $i\in \kappa$ then the resulting
chains are not adjacent and the edge must be removed.
Figure \ref{F:type2} shows the result
of combining $C_1$ with both $C_2$ and $C_2'$ using the shuffles associated with
3-sequences $(0,2,2)$ and $(0,1,1)$.  Chains $C_2$ and $C_2'$ differ at rank 1,
but $(0,2,2)$ does not contain an element 1, so the shuffle does not stretch the
diamond and the resulting chains are adjacent in $L_1 \times L_2$.
However, the shuffle associated with $(0,1,1)$ stretches the diamond by inserting
two edges from $C_1$, so this Type 2 edge must be removed from
$\widetilde{\Gamma}(L_1 \times L_2)$.

\textbf{Type 3 edges}.  $C_1 \sim C_1'$ in $L_1$, $C_2=C_2'$,
and $\kappa = \kappa'.$
As in the previous case, we must first identify the rank $i$
where $C_1$ and $C_1'$ differ.
If $i\in \lambda$
then the chains are not adjacent in
$L_1 \times L_2$ and the edge is removed from $\widetilde{\Gamma}(L_1 \times L_2)$.

This completes the determination of which edges to remove from
$\widetilde{\Gamma}(L_1 \times L_2)$, resulting in the final
graph $\Gamma(L_1 \times L_2)$.

%%%%%%%%%%%%%%%%%%%%%%%%%%%%%%%%%%%%%%%%%%%%%%%%%%%%%%%%%
\section{The Boolean Lattice}
\label{S:Graph}

\begin{figure}[ht]
\begin{center}
\includegraphics[height=3.5in]{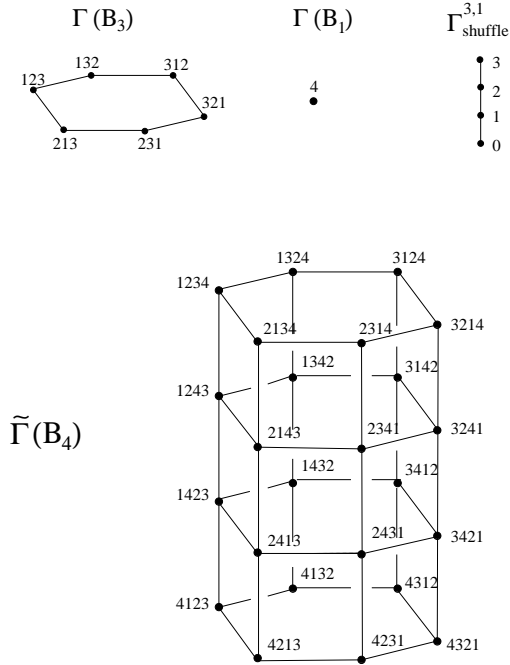}
\caption{The intermediate graph $\widetilde{\Gamma}(B_4)$.}
\label{F:GammaTilde}
\end{center}
\end{figure}

While in \cite{mainpaper}, it was shown that for any connected graphs
$\Gamma_1$ and $\Gamma_2$, we have that $A^G_1(\Gamma_1 \Box \Gamma_2)
\cong A^G_1(\Gamma_1) \times A^G_1(\Gamma_2),$
it is not immediately clear from the construction of $\Gamma(L_1 \times L_2)$
if there is an easily defined relationship between the groups
$A^G_1(\Gamma(L_1))$, $A^G_1(\Gamma(L_2))$, and
$A^G_1(\Gamma(L_1 \times L_2))$.  However, applying the construction defined in
the previous section to $B_n$ leads to a better understanding of $\Gamma(B_n)$
that allows us to use combinatorial methods to compute the abelianization of
its $A_1$-group.  We now reconstruct $\Gamma (B_n)$ and characterize all of its
primitive $6$-cycles.

In Figure \ref{F:GammaTilde}, the vertices of $\Gamma(B_3)$ and $\widetilde{\Gamma}(B_4)$ are
labeled with permutations written in one line notation.  The graph $\Gamma(B_1)$
consists of a single vertex, thus there are no Type 2 edges in
$\widetilde{\Gamma}(B_n)$ and we can simply consider the 1-sequences of
$\Gamma^{3,1}_{shuffle}$.  Labeling $\Gamma(B_1)$ with the element
4 enables us to see the relationship between the 1-sequences labeling the shuffle
graph and the position of 4 in the resulting permutations in $S_4$.

\begin{figure}[ht]
\begin{center}
\includegraphics[height=3.2in]{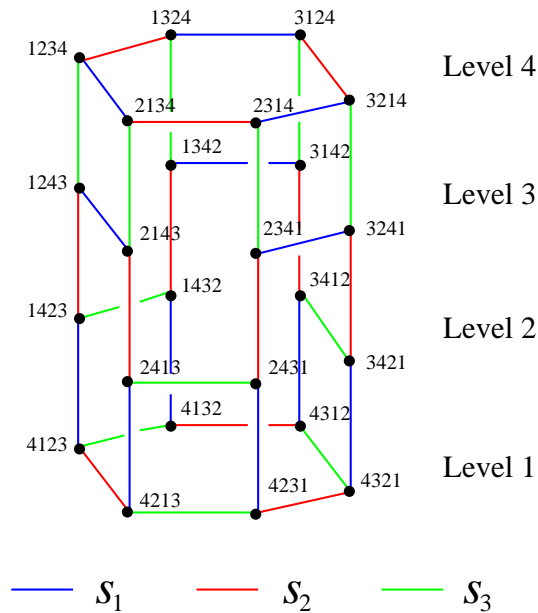}
\caption{The final graph $\Gamma(B_4)$.}
\label{F:GammaB4}
\end{center}
\end{figure}

The graph in Figure \ref{F:GammaB4} is another representation of the
permutahedron we saw in Figure \ref{F:permutahedron}. The vertices are
labeled with permutations of $S_4$, written in one line notation, and
each edge corresponds to a (simple) transposition.
$\Gamma(B_4)$ illustrates the following (most of them well-known) properties
of $\Gamma(B_n).$ All properties easily follow from the definition of the graph.

\begin{center}
\large \textbf{Properties of $\mathbf{\Gamma(B_n)}$}
\end{center}

\begin{enumerate}

\item
The graph $\Gamma(B_n)$ is (\textit{n}\,-1)-regular, with each vertex
incident to precisely one edge for each of the $n-1$ transpositions
$s_i\in \bf{S},$ $1\leq i \leq n-1.$ Label each edge with its associated transposition.

\item
Let {\it level} $i$ be the set of all permutations $\sigma \in S_n$
such that $\sigma^{-1}(n)=i.$ The graph has $n$ levels,
and each level was initially a copy of $\Gamma(B_{n-1})$ before we removed edges
from $\widetilde{\Gamma}(B_n)$.

\item
The graph $\Gamma(B_n)$ is bipartite, with vertices partitioned into even and odd
permutations, and all cycles in the graph are of even length.

\item
The sequence of transpositions labeling the edges of a cycle in
$\Gamma(B_n)$ form a representation of the identity in $S_n$.
Each 4-cycle in the graph corresponds to
$(s_ks_j)^2$ for some $1 \leq k,\ j \leq n-1,$ where $|j-k|\geq 2.$
Each primitive 6-cycle corresponds to
$(s_js_{j+1})^3$ for some $1\leq j \leq n-2$.

\end{enumerate}

Due to the structure of $S_n$, all other cycles of length $\geq 8$
can be expressed as the concatenation of 4-
and 6-cycles; we can therefore limit our investigation to primitive 6-cycles.
We want to count the $G$-homotopy equivalence classes of 6-cycles in $B_n$,
which yields the rank of $A^G_1(\Gamma(B_n))^{ab}$.  The following definition gives us
an additional description of edges and cycles in $\Gamma(B_n)$ that will help us
determine equivalence classes of 6-cycles.

\begin{definition}
$ $
\begin{enumerate}

\item
An edge in $\Gamma(B_n)$ between $\sigma$ and $\tau$ is \textit{horizontal} if
$\sigma^{-1}(n)=\tau^{-1}(n),$ or \textit{vertical} if $\sigma^{-1}(n)=\tau^{-1}(n)\pm 1.$

\item
All vertices in a \textit{horizontal} 6-cycle are at the same level.
A \textit{vertical} 6-cycle contains two vertices in each of three consecutive levels.

\end{enumerate}
\end{definition}

We note that all vertical edges between levels $i$ and $i+1$ are labeled
with $s_i.$  We identify each vertical 6-cycle with the middle
of its three levels.  For example, 1243-2143-2413-4213-4123-1423 is a vertical
6-cycle at level 2 in $\Gamma(B_4)$, and its edges are labeled with $s_1$ and $s_2$.

Let $C$ be a cycle in $\Gamma(B_{n-1})$.  There are $n$ copies of $C$ in
$\widetilde{\Gamma}(B_n)$, one in each level of the graph.  Each copy of $C$
in $\widetilde{\Gamma}(B_n)$
is a horizontal cycle that is connected to the copies in neighboring levels by
vertical edges, forming a net of vertical 4-cycles connecting all of the copies.
For example, we see in Figure \ref{F:GammaTilde} that the 6-cycles forming
levels 1 and 4 are connected by such a net.  Removing vertical edges from
$\widetilde{\Gamma}(B_n)$ to form $\Gamma(B_n)$ may remove edges from this net,
however, it is easy to see that the copies of $C$ remain connected by a net of
vertical 4- and 6-cycles, as seen in Figure \ref{F:GammaB4}.

%%%%%%%%%%%%%%%%%%%%%%%%%%%%%%%%%%%%%%%%%%%%%%%%%%%%%%%%%%%%%%%%%%%%%%%%%%%%%%%%%%%%%%
\section{Equivalence Classes}
\label{S:Classes}

In this section, we describe how to distinguish between different $G$-homotopy equivalence classes in $\Gamma(B_n)$ so that we may count them.  Specifically, we consider graph maps whose images in $\Gamma(B_n)$ are primitive 6-cycles connected to the base vertex by a path, and we determine when they are $G$-homotopic to one another.  We denote this relationship by $C_1 \simeq_G C_2$, referring to the 6-cycles in the images rather than the graph maps themselves.

First, we show that if two 6-cycles are in the same equivalence class then they
are associated with the same pair of transpositions, $s_i$ and $s_{i+1}$, for
some $i$, $1\leq i \leq n-1$.  We then prove a stronger theorem:  6-cycles $C_1$
and $C_2$ are in the same equivalence class if and only if they differ by a
sequence of transpositions disjoint from $s_i$ and $s_{i+1}$.
This theorem, when combined with our
new understanding of the structure of $\Gamma(B_n)$, gives us the means to
describe the equivalence classes, and to find a formula for the rank of
$A_1^G(\Gamma(B_n))^{ab}$.

\begin{figure}[ht]
\begin{center}
\includegraphics[height=1.2in]{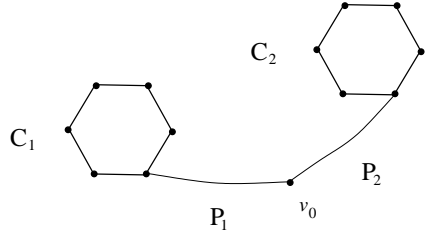}
\caption{Two 6-cycles based at $v_0$.}
\label{F:BasedLoops}
\end{center}
\end{figure}

Let $C_1$ and $C_2$ be distinct primitive 6-cycles in $\Gamma(B_n)$.
Let $\sigma_0$ be the base permutation and let $P_1$ and $P_2$ be paths
from $\sigma_0$ to $C_1$ and $C_2,$ respectively. Then
by definition  $C_1\simeq_G C_2$
if and only if there exists a $G$-homotopy grid such that
the first row is $P_1 C_1 P_1^{-1}$ and the last row is $P_2 C_2 P_2^{-1}$.

Suppose that $C_1 \simeq_G C_2$ and we have a $G$-homotopy grid from
$P_1 C_1 P_1^{-1}$ to $P_2 C_2 P_2^{-1}$.  We label each vertex in the
grid with the corresponding permutation in $S_n$.  Recall that a $G$-homotopy
is itself a graph map that preserves adjacency, thus two vertices, $\sigma_1$ and $\sigma_2$,
are adjacent in the grid if and only if
$\sigma_1 = \sigma_2$, or $\sigma_2=\sigma_1 s_i$
for some $s_i$, $1 \leq i \leq n-1$.  If $\sigma_2=\sigma_1 s_i$,
then we label the edge with $s_i$; if $\sigma_1=\sigma_2$ the edge remains unlabeled.
Each row in the grid is a loop based at $\sigma_0$, to which we associate a word
in $S_n$ formed
by the sequence of transpositions labeling the edges from left to right in the row.
Let

$$\phantom{and} \; \; W_1=s_{i_1} s_{i_2} \dots s_{i_m} (s_i s_{i+1})^3 s_{i_m} \dots s_{i_2} s_{i_1}$$
$$\textnormal{and} \; \; W_2=s'_{i_1} s'_{i_2} \dots s'_{i_n} (s_j s_{j+1})^3 s'_{i_n} \dots s'_{i_2} s'_{i_1}$$

\noindent
be the words associated with $P_1 C_1 P_1^{-1}$ and $P_2 C_2 P_2^{-1}$, respectively.
We note that each word is a non-reduced representation of the identity
and $\sigma_0 W_1$ and $\sigma_0 W_2$ are \textit{based words}.
Thus a $G$-homotopy between primitive 6-cycles $C_1$ and $C_2$ corresponds to
a sequence of operations that transform $W_1$ into $W_2$.

We now consider the possible changes that we may make from one row to the next
which preserve the $G$-homotopy relation, and describe each change in terms
of operations on the associated words.  Two graph maps are $G$-homotopic if
they differ by 3- and 4-cycles, however, $\Gamma(B_n)$ contains no 3-cycles.
Furthermore, we saw in Figure \ref{F:homotopyFtoG} that it is not possible
to contract a 4-cycle in a single step, thus the changes described below are the
only permissible changes.

\begin{enumerate}

\item[({\bf T1})]  \textbf{Repeating vertices.}
The $G$-homotopy relation is preserved by the repetition of a vertex in the image of a row, and the corresponding word does not change.

\begin{figure}[h]
\begin{center}
\includegraphics[height=1.3in]{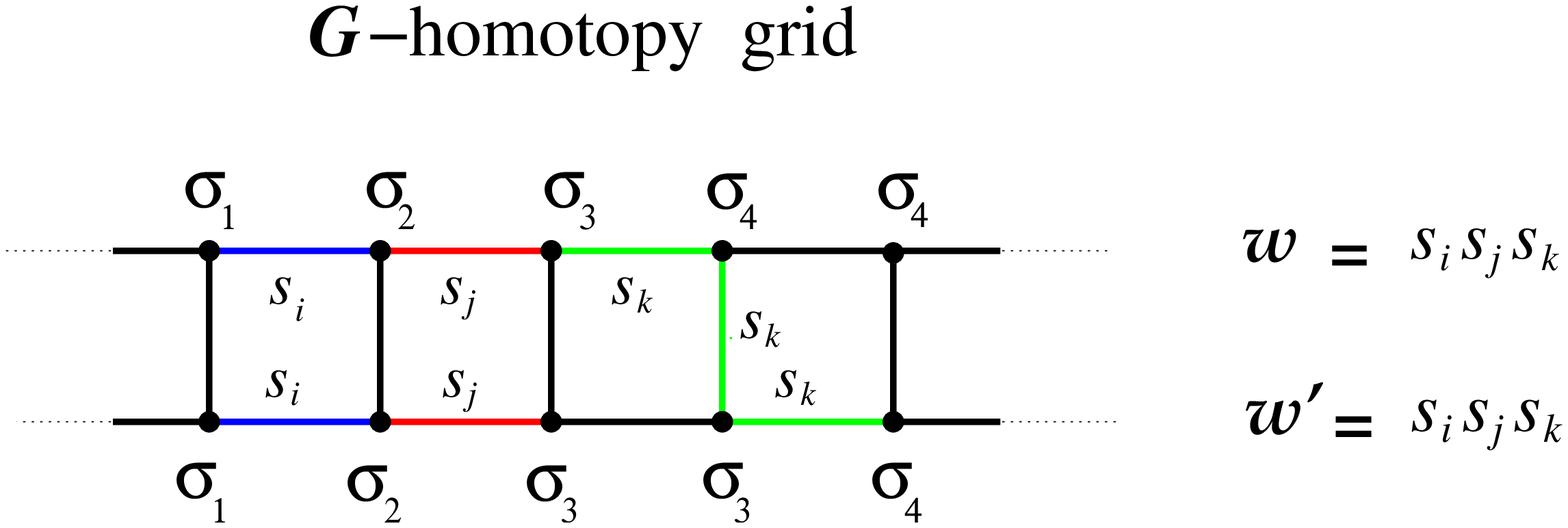}
\caption{Repeating a vertex.}
\label{F:RepeatVertex}
\end{center}
\end{figure}

\item[({\bf T2})] \textbf{Traversing an edge in both directions.}
Traversing an edge once in each direction also preserves the $G$-homotopy relation and is equivalent to inserting $s_i^2$ into the associated word.  Similarly, we may reverse this by removing such an edge, and deleting $s_i^2$ from the word.

\begin{figure}[h]
\begin{center}
\includegraphics[height=2.6in]{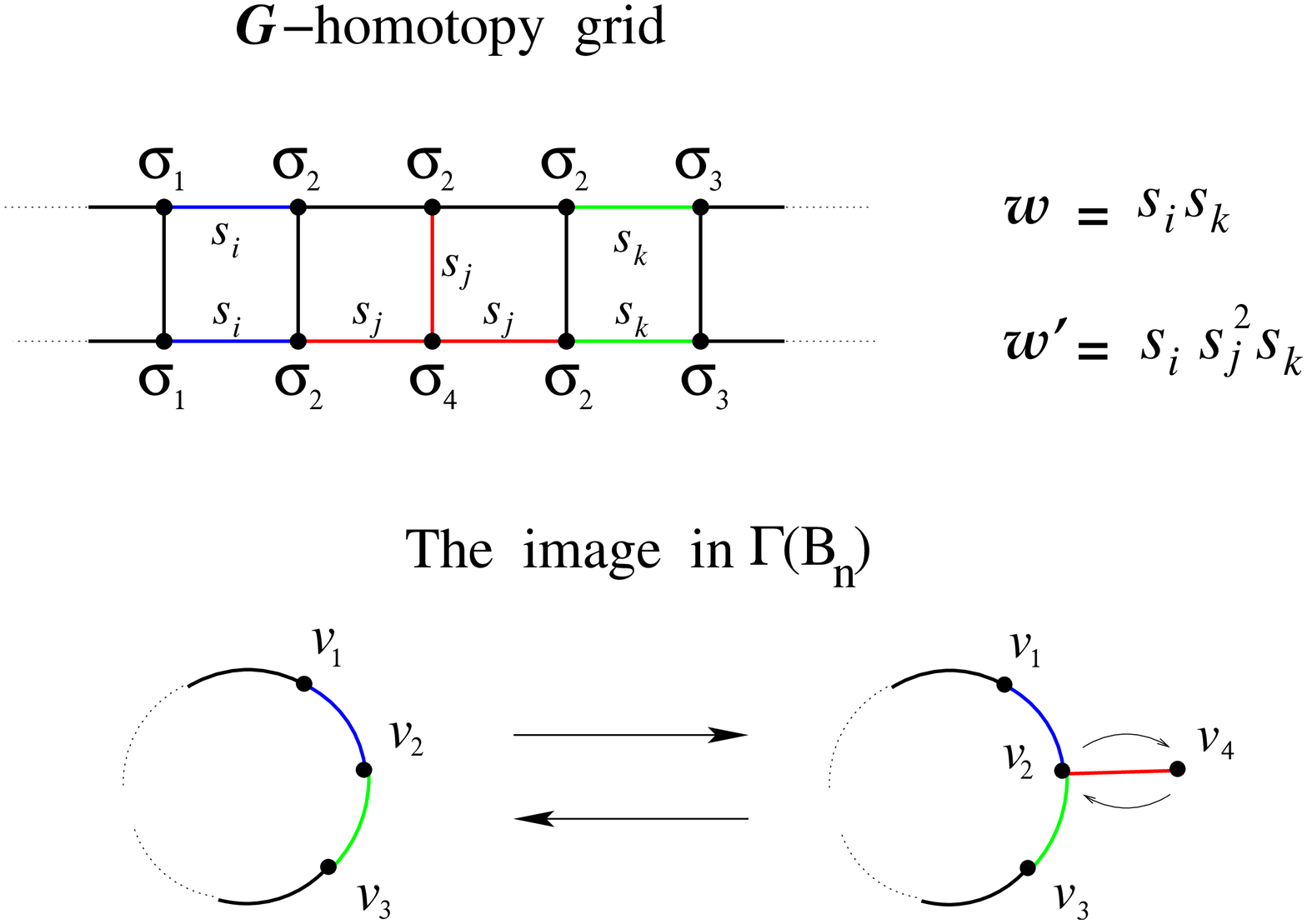}
\caption{Traversing an edge in both directions.}
\label{F:BothDirections}
\end{center}
\end{figure}

\item[({\bf T3})] \textbf{Exchanging pairs of edges of a 4-cycle.}
Since 4-cycles are $G$-homotopic to the identity, we can replace two
consecutive edges of a 4-cycle with the other two edges in the cycle.
Recall that a 4-cycle is associated with $s_j$ and $s_k$, with $|j-k| \geq 2$.
We see that the effect of this change is to commute $s_j$ and $s_k$ in the word.
See Figure \ref{F:Exchange2-2}.
\begin{figure}[h]
\begin{center}
\includegraphics[height=2.6in]{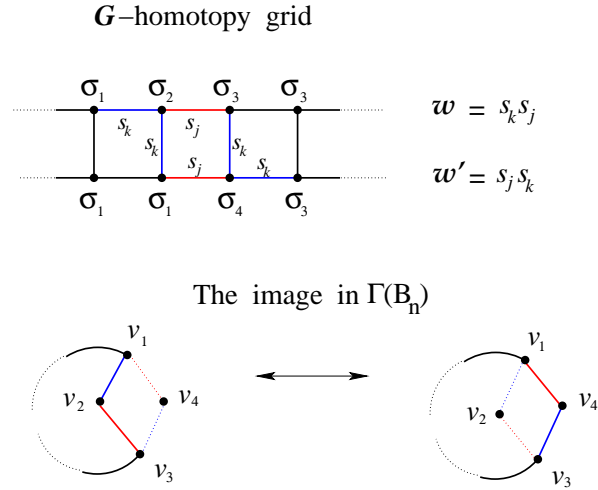}
\caption{Exchanging pairs of edges.}
\label{F:Exchange2-2}
\end{center}
\end{figure}

\end{enumerate}

A key feature of each of the changes described above is that they
preserve the parity of the number of transpositions $s_i$ (for each $1\leq i\leq n-1$) in
the associated words.  We also note that these changes do not include the
use of the relation $(s_j s_{j+1})^3=1$ because that would involve exchanging
consecutive edges of a primitive 6-cycle with the remaining edges.  This would
obviously not preserve $G$-homotopy because any cycle of length $\geq5$ is not
contractible to a single vertex.

Note that the above two types of operations, $(\textbf{T2})$ and $(\textbf{T3})$
generate an equivalence relation on the set $\sigma_0 \textbf{W}$
of (based) loop-words, where
$$\textbf{W}=\{W=s_{j_1} s_{j_2} \ldots s_{j_{2k}}, \forall k \geq 1 \:| \;
s_{j_i}\in {\bf{S}}, \ W=1\}.$$ Two based words $\sigma_0 W$ and $\sigma_0 W'$
are equivalent, $ \sigma_0 W \sim \sigma_0 W',$ if one can be obtained from the
other by a sequence of $(\textbf{T2})$ and $(\textbf{T3})$ operations.
Furthermore, two 6-cycles based at $\sigma_0$, $C_1$ and $C_2$, are $G$-homotopic if
and only if their corresponding words,
$\sigma_0 W_1$ and $\sigma_0 W_2$  are equivalent.
Therefore we can continue our investigation of equivalence classes of 6-cycles
using the language of graphs or words interchangeably.

\begin{proposition}
\label{P:SameTranspositions}
Let $C_1$ and $C_2$ be two $G$-homotopic  primitive 6-cycles in $\Gamma(B_n)$.
If $C_1 \simeq_G C_2$, then $C_1=C_2=(s_is_{i+1})^3,$
for some $i$, $1 \leq i \leq n-2$.
\end{proposition}

\begin{proof}
For $i=1,2,$ let $P_i$ be a path from $\sigma_0$ to $C_i.$
By  assumption $P_1C_1P^{-1}_1 \simeq_G P_2C_2P^{-1}_2,$ and the corresponding words
$w_1(s_is_{i+1})^3w^{-1}_1,$ $w_2(s_js_{j+1})^3w^{-1}_2$ are equivalent.
This means that we must be able to transform the first word into the second
using only $(\textbf{T2})$ and $(\textbf{T3})$ operations.
A simple parity argument ($s_i$ appears an odd number of times in
$w_1(s_is_{i+1})^3w^{-1}_1$) shows that this is possible only if $i=j.$
\end{proof}

Association with the same pair of transpositions is a necessary condition
for
\linebreak
$G$-homotopy of 6-cycles, but it turns out not to be sufficient.
For $1\leq i \leq 6,$ let $\sigma_i$ and $\gamma_i$ be the vertices of
the cycles $C_1$ and $C_2$ respectively. Since $\Gamma(B_n)$
is connected there exists a path from $C_1$ to $C_2.$
Let $\tau_i\in \bf{S}$ and let $\tau_1, \tau_2, \ldots \tau_k$
be a shortest path from $C_1$ to $C_2.$
Moreover, assume that $\gamma_1=\sigma_1\tau_1\tau_2 \ldots \tau_k,$
and that we go around the cycles $C_1$ and $C_2$
in the $s_i$ direction.
That is, $\sigma_2=\sigma_1s_i$, $\sigma_3=\sigma_1s_is_{i+1}$ and so on.
Similarly $\gamma_2=\gamma_1s_i,$ $\gamma_3=\gamma_1s_is_{i+1},$ etc.
With this notation in mind we have
the following theorem.

\begin{theorem}
\label{T:DifferByTranspositions}
Let $C_1$ and $C_2$ be two distinct primitive 6-cycles in $\Gamma(B_n)$.
Then $C_1 \simeq_G C_2$ iff there exists $k \geq 1,$ and transpositions
$\tau_1\tau_2...\tau_k$
disjoint from $s_i$ and $s_{i+1},$
such that $\gamma_i=\sigma_i \tau_1 \dots \tau_k$, for all $i=1,\ldots 6.$

\end{theorem}

\begin{proof}The first part of the proof is constructive:
assuming $C_2=C_1 \tau_1 \dots \tau_k$, we are able to construct a $G$-homotopy
that connects $C_1$ to $C_2$ by a net of 4-cycles of type $(s_i\tau_j)^2$ or
$(s_{i+1}\tau_j)^2.$
Note that if $\tau_j$ is not disjoint from $s_i$ and $s_{i+1},$ we do
not get $4$-cycles.
Figure \ref{F:6Cycles} is the image of
one such $G$-homotopy from $C_1$ to $C_2=C_1\tau_1\tau_2\tau_3$.

\vskip .2in

\begin{figure}[ht]
\begin{center}
\includegraphics[height=2.4in]{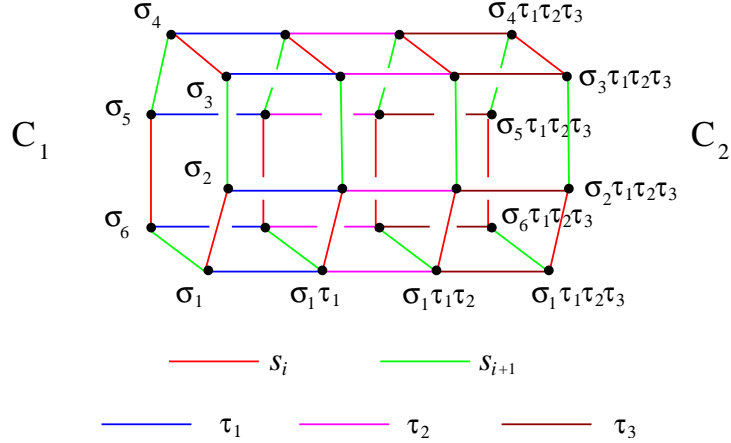}
\caption{A $G$-homotopy from $C_1$ to $C_2$.}
\label{F:6Cycles}
\end{center}
\end{figure}

For the second part of the proof assume that $C_1$ and $C_2$ are $G$-homotopic
primitive 6-cycles, thus $C_1=C_2=(s_is_{i+1})^3,$
for some $1 \leq i \leq n-1$.
Let $P$ be a shortest path from $C_1$ to $C_2$,
with the vertices of both cycles as described
above the theorem.
See Figure \ref{F:ShortestPath}.

While neither $C_1$ nor $C_2$ can be contracted, the loop $C_1 P C_2^{-1} P^{-1}$
is contractible to a single vertex.  Let $w=\tau_1 \tau_2 \dots \tau_k$
be the word corresponding to $P$,
thus $$W=(s_i s_{i+1})^3 w (s_{i+1} s_i)^3 w^{-1}$$ corresponds to
$C_1 P C_2^{-1} P^{-1}$. Since as a permutation $W=1$
we must be able to reduce the word $W$ to the empty word
using only the relations
$s_j^2=1$ and $s_j s_k = s_k s_j$ if $|j-k| \geq 2.$

\begin{figure}[ht]
\begin{center}
\includegraphics[height=1.6in]{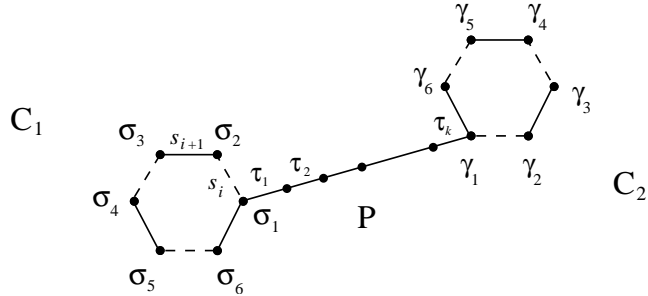}
\caption{The loop $C_1 P C_2^{-1} P^{-1}$ can be contracted to a single vertex.}
\label{F:ShortestPath}
\end{center}
\end{figure}

Our goal is to show that $w$ consists solely of transpositions disjoint from
$s_i$ and $s_{i+1}$, and consequently that $\gamma_i=\sigma_i \tau_1 \tau_2 \dots \tau_k,$
for all $1\leq i \leq 6.$
This is proven by first showing that we must be able to
reduce $W$ to the empty word without having to insert $s_j$, for any $1\leq j\leq n-1.$
Next,
we consider the types of transpositions that may occur in $w$, and show
that it may only contain those that are disjoint from $s_i$ and $s_{i+1}.$
This part of the proof requires checking many cases.
Once we have shown that $w=\tau_1 \tau_2 \dots \tau_k$, where each of the $\tau_j$ are disjoint from $s_i$ and $s_{i+1}$, then we can see that

\begin{displaymath}
    \begin{array}{l}
    \gamma_1 = \sigma_1 \tau_1 \tau_2 \dots \tau_k \textnormal \quad \textnormal{and}\\

    \gamma_2 = \gamma_1 s_i\\
    \phantom {\gamma_1} = \sigma_1 \tau_1 \tau_2 \dots \tau_k  s_i\\
    \phantom {\gamma_1} = \sigma_1 s_i \tau_1 \tau_2 \dots \tau_k\\
    \phantom {\gamma_1} = \sigma_2 \tau_1 \tau_2 \dots \tau_k.
  \end{array}
\end{displaymath}

By a similar argument, $\gamma_j    = \sigma_j \tau_1 \tau_2 \dots \tau_k$ for
$3 \leq j \leq 6.$

\begin{center}
\textbf{ Insertion of $\mathbf{s_j^2}$ is not needed.}
\end{center}

Suppose we insert $s_j^2$, at some step in the process.
By the end of the process, each of these two $s_j$  will have been removed.
If those two $s_j$ were removed
together as a single pair, then they did not assist us in removing other
occurrences of $s_j$, so it was not necessary to insert them at all.

If the $s_j$ were removed separately, each paired with another occurrence of $s_j$,
then whatever commuting that was done to put them into position next to their
new partners could have been done in the opposite direction by the partner
transpositions, which we could then remove.  Thus we must be able to
reduce $W$ to the empty word without inserting transpositions.

\begin{center}
\textbf{The word \textit{w} does not contain $\mathbf{s_{i-1}}$ or $\mathbf{s_{i+2}}$.}
\end{center}

Without loss of generality, suppose there is at least one $s_{i-1}$ in $w.$ Since $w$ is a shortest path
the only way to cancel that $s_{i-1}$ is by pairing it with the first $s_{i-1}$ in $w^{-1}.$ But there is an odd number of
transpositions $s_i$ between the last occurrence of $s_{i-1}$ in $w$ and
its first occurrence in $w^{-1}$; three from $(s_{i+1} s_1)^3$ and an even number,
if any, from $w$ and $w^{-1}$.  Thus, even if there were to be some cancelation
there will be at least one occurrence of $s_i$ left between the $s_{i-1}$,
preventing their pairing.

\begin{center}
\textbf{The word \textit{w} does not contain $\mathbf{s_i}$ or $\mathbf{s_{i+1}}$.}
\end{center}

Assume now that there are some $s_i$ or $s_{i+1}$ in $w.$
Consider the subsequence of all the transpositions
$s_i$ and $s_{i+1}$ in $w$.  We show that if this subsequence is not empty
then $P$ is not a shortest path.

First, we show that a subsequence cannot have consecutive occurrences of either
$s_i$ or $s_{i+1}$.  Next, we show that the subsequences must be of length shorter than
six.  We then show that transforming $W$ to the empty
word is not possible if the alternating subsequence is of odd length.
Finally, we deal with the subsequences of even length.

\textbf{Subsequence of $\mathbf{s_i}$, $\mathbf{s_{i+1}}$ alternates.} Suppose that
the subsequence in $w$ consisting of all transpositions $s_i$ and $s_{i+1}$
contains a consecutive pair of the same transposition, say $s_i.$
Since $s_i$
commutes with all transpositions in $w$ different than
$s_{i+1}$, we can commute those two $s_i$ within $w$ until they are adjacent at which
point they can be removed. But this yields a shorter word $w$, and thus a shorter
path $P.$

\textbf{Length of subsequence is less than six.} Next we consider which lengths
of alternating subsequences of $s_i$ and $s_{i+1}$ are possible in $w$.
If we have an alternating subsequence of length six or longer, then we
can commute $s_i$ and $s_{i+1}$ with the other transpositions in $w$ until
we have the subword $(s_is_{i+1})^3.$ The
new word corresponds to a path from $C_1$ to $C_2$ that contains a 6-cycle.
The same path without the 6-cycle would be shorter.

\textbf{Odd subsequences.} For subsequences of odd length,
Figure \ref{F:NewCancellation} illustrates the possible subsequences consisting
of $s_i$ and $s_{i+1}$ in $W$.  This includes the six transpositions from
each of $C_1$ and $C_2$, as well as those from $w$ and $w^{-1}$.  In each case,
we see that we are not
able to remove all transpositions in the entire subsequence, contradicting the
assumption that $C_1$ and $C_2$ are $G$-homotopic.

\begin{figure}[ht]
\begin{center}
\includegraphics[height=2.1in]{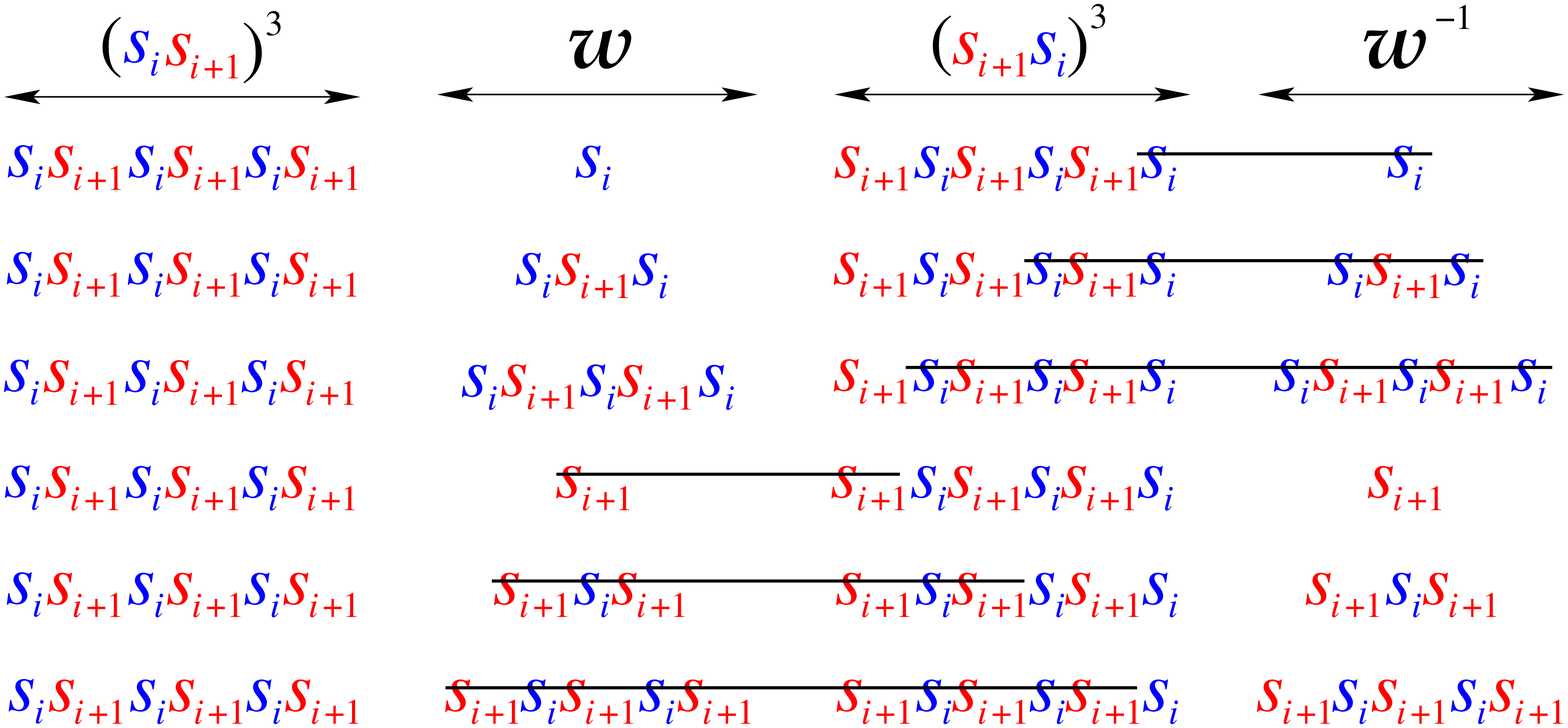}
\caption{We cannot reduce $W=(s_i s_{i+1})^3 w (s_{i+1} s_i)^3 w^{-1}$ to the empty word when subsequences are of odd length.}
\label{F:NewCancellation}
\end{center}
\end{figure}

\textbf{Even subsequences.} Assume that the alternating subsequence in $P$
is of length two or four.  To illustrate this, consider the case where
$$w=\tau_1 \tau_2 s_i \tau_3 s_{i+1} \tau_4.$$  We can commute these transpositions within $w$ to obtain
$$w'=s_i s_{i+1}  \tau_1  \tau_2 \tau_3 \tau_4,$$ with
the transpositions $\tau_j$ in the same relative order as they were in $w$.
Note that as permutations in $S_n,$ $w=w'.$

Recall that $\sigma_j$ and $\gamma_j$, $1 \leq j \leq 6$, are the
permutations associated with the vertices in $C_1$ and $C_2$,
respectively, with $\sigma_1$ the first vertex in $P$,
and $\gamma_1$ the last.  Our reordered word $w'$ also
corresponds to a path from $\sigma_1$ to $\gamma_1$ in $\Gamma(B_n)$,
thus we may write

\begin{displaymath}
    \begin{array}{ll}
    \gamma_1 = \sigma_1 s_i s_{i+1}  \tau_1  \tau_2 \tau_3 \tau_4\\
    \phantom {\gamma_1} = \sigma_3 \tau_1  \tau_2 \tau_3 \tau_4.
    \end{array}
\end{displaymath}

\noindent and similarly we have
\begin{displaymath}
    \begin{array}{ll}
    \gamma_2 = \sigma_4 \tau_1  \tau_2 \tau_3 \tau_4\\
    \gamma_3 = \sigma_5 \tau_1  \tau_2 \tau_3 \tau_4\\
    \gamma_4 = \sigma_6 \tau_1  \tau_2 \tau_3 \tau_4\\
    \gamma_5 = \sigma_1 \tau_1  \tau_2 \tau_3 \tau_4\\
    \gamma_6 = \sigma_2 \tau_1  \tau_2 \tau_3 \tau_4.
    \end{array}
\end{displaymath}

\begin{figure}[ht]
\begin{center}
\includegraphics[height=2.5in]{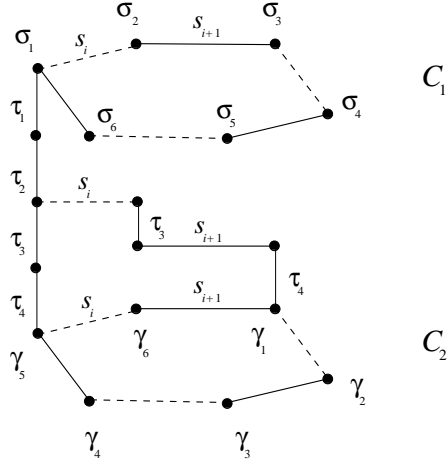}
\caption{A shorter path without $s_i$ and $s_{i+1}.$}
\label{F:NewPath}
\end{center}
\end{figure}

\noindent We see in Figure \ref{F:NewPath} that the word $\tau_1  \tau_2 \tau_3 \tau_4$ corresponds to a shorter path from $\sigma_1$ to $C_2.$ This new path now meets $C_2$ at $\gamma_5$, but it still satisfies the condition that the first edge traversed in $C_2$
is $s_i$, thus contradicting the assumption that $P$ is a shortest path.
A similar argument holds true for each of the remaining
three cases.
\end{proof}

The proof of Theorem \ref{T:DifferByTranspositions} relies on the definition
of a $G$-homotopy from $C_1$ to $C_2$, and the limitations on the types of
changes we are able to make to the words associated with each row of the
$G$-homotopy grid.  We can combine this theorem with the properties
of $\Gamma(B_n)$ to obtain the following corollaries about
equivalence classes of primitive 6-cycles.

\begin{corollary}
Horizontal and vertical 6-cycles are in different equivalence classes.
\end{corollary}

\begin{proof}
Suppose $C_1$ is a horizontal 6-cycle at level $i$. As observed earlier,
for all vertices $\sigma$ in $C_1$, $\sigma^{-1}(n)=i$.
Furthermore, if $C_2=C_1\tau_1 \dots \tau_k$, then
$(\sigma \tau_1 \dots \tau_k)^{-1}(n)=\tau_k \dots \tau_1 (i)$ and
consequently $C_2$ is also a horizontal 6-cycle.  Similarly, if $C_1$
is a vertical 6-cycle, then $C_2$ is a vertical 6-cycle.  Therefore we may
count horizontal and vertical equivalence classes separately.
\end{proof}

\begin{corollary}
Vertical 6-cycles at different levels of $\Gamma(B_n)$ are in different equivalence classes.
\end{corollary}

\begin{proof}
This is a direct consequence of Proposition \ref{P:SameTranspositions} and the observation we made in Section \ref{S:Graph} that a vertical 6-cycle at level $i$, $2 \leq i \leq n-1$ is associated with transpositions $s_{i-1}$ and $s_i$.
\end{proof}

\begin{corollary}
For all $2 \leq i \leq n-1,$ the number of vertical equivalence classes
at level $i$ in $\Gamma(B_n)$ is given by
$$\frac{(n-1)!}{2(i-2)!(n-i-1)!}=\binom{n-1}{i} \binom{i}{2}.$$
\end{corollary}

\begin{proof}
There are $(n-1)!$ vertices in each level of $\Gamma(B_n)$, and each vertical 6-cycle at level $i$ contains 2 vertices in that level, thus there are $\frac{(n-1)!}{2}$ vertical 6-cycles at level $i$.  The number of 6-cycles in each vertical equivalence class at level $i$ is equal to the order of the subgroup of $S_n$ generated by transpositions disjoint from $s_{i-1}$ and $s_i$, which is $(i-2)!(n-i-1)!$, thus the number of equivalence classes at level $i$ is $\frac{(n-1)!}{2(i-2)!(n-i-1)!}$.
\end{proof}

Consequently, the total number of vertical equivalence classes is

\begin{eqnarray*}
\sum_{i=2}^{n-1}\binom{n-1}{i} \binom{i}{2} &=& \sum_{i=2}^{n-1} \frac{(n-1)!}{2(i-2)!(n-i-1)!}\\
\phantom{\sum_{i=2}^{n-1}\binom{n-1}{i} \binom{i}{2}} &=& \sum_{i=0}^{n-3} \frac{(n-1)!}{2(i)!(n-i-3)!}\\
\phantom{\sum_{i=2}^{n-1}\binom{n-1}{i} \binom{i}{2}} &= &\left( \sum_{i=0}^{n-3}\ \binom{n-3}{i}\right) \frac{(n-1)(n-2)}{2} \\
\phantom{\sum_{i=2}^{n-1}\binom{n-1}{i} \binom{i}{2}} &=& 2^{n-3}\binom{n-1}{2}.
\end{eqnarray*}

\begin{corollary}
The number of horizontal equivalence classes in $A^G_1(\Gamma(B_n))^{ab}$ is equal to the rank of $A^G_1(\Gamma(B_{n-1}))^{ab}$.
\end{corollary}

\begin{proof}
We count the horizontal equivalence classes in level $n$, which remains a copy of
$\Gamma(B_{n-1})$ even after we removed edges from $\widetilde{\Gamma}(B_n)$
in the construction of $\Gamma(B_n)$.  Let $C$ be a horizontal 6-cycle in the graph
that is not at level $n$.  Recall that there is a copy, $C'$, of $C$ at level $n$
of $\Gamma(B_n)$, and there is a net of vertical 4- and 6-cycles connecting $C$
to $C'$.  While $C$ may not be in the same equivalence class as $C'$, its equivalence class
can be expressed in terms of those of $C'$ and the
vertical 6-cycles in the net.
\end{proof}

We are now able to see why $A^G_1(\Gamma(B_4))^{ab}$ has seven generators
rather than eight.  In $\Gamma(B_4)$, the horizontal 6-cycle at level 4 and
the six vertical 6-cycles are each in distinct equivalence classes; the
remaining 6-cycle equivalence class at level 1 is clearly expressible in terms of the
previous classes.

\begin{theorem}
For all $n\geq 1,$
$$\text{rank\ } A^G_1(\Gamma(B_n))^{ab}= 2^{n-3}(n^2-5n+8)-1.$$
\end{theorem}

\begin{proof}
From Theorem $5.2$ and all of its corollaries we have that there are
$\sum_{k=1}^{n}2^{k-3} \binom{k-1}{2}$ equivalence classes of 6-cycles in
$\Gamma(B_n)$.  Any other cycle in $\Gamma(B_n)$ of length $\geq 8$ can
be expressed as the concatenation of 4-cycles and 6-cycles, so those $6$-cycles classes
generate the free group $A_1^{n-3}(\Delta(B_n))^{ab}$.

Let $f(n) = \sum_{k=1}^{n}2^{k-3} \binom{k-1}{2}= \text{rank\ }A^G_1(\Gamma(B_n))^{ab}$. Then

\begin{eqnarray*}
\sum_{n\geq 0} f(n)x^n&=& \sum_{n \geq 0}  \sum_{k=0}^{n} 2^{n-k-3}\binom{n-k-1}{2} x^n\\
%\phantom{F(x)}&=&\sum_{k\geq 0} \sum_{n} 2^{n-k-3}\binom{n-k-1}{2} x^n\\
\phantom{\sum_{n\geq 0} f(n)x^n}&=&\sum_{k\geq 0}2^{-2}x^{k+1}\sum_{n} 2^{n-k-1}\binom{n-k-1}{2} x^{n-k-1}\\
\phantom{\sum_{n\geq 0} f(n)x^n}&=&\sum_{k\geq 0}\frac{1}{4}x^{k+1}
\sum_{r} \binom {r}{2}(2x)^r\\
\phantom{\sum_{n\geq 0} f(n)x^n}&=&\sum_{k\geq 0}\frac{1}{4}x^{k+1}\frac {(2x)^2}{(1-2x)^3}\\
%\phantom{\sum_{n\geq 0} f(n)x^n}&=&\sum_{k\geq 0}\frac{x^{k+3}}{(1-2x)^3}\\
\phantom{\sum_{n\geq 0} f(n)x^n}&=&\frac{x^3}{(1-2x)^3}\sum_{k\geq 0}x^k\\
\phantom{\sum_{n\geq 0} f(n)x^n}&=&\frac{x^3}{(1-2x)^3(1-x)}.
\end{eqnarray*}

Using partial fraction decomposition, we
see that $$f(n)=2^{n-3}(n^2-5n+8)-1.$$
\end{proof}

%%%%%%%%%%%%%%%%%%%%%%%%%%%%%%%%%%%%%%%%%%%%%%%%%%%%%%%%%%%%%%%%%%%%%%%%%%%%%%%%%%%%%%

In 2001, Babson observed that
attaching 2-cells to the 4-cycles in $\Gamma(B_n)$ results in a topological space
that is homotopy equivalent to the complement (in $\mathbb{R}^n$) of the
3-equal hyperplane arrangement.  But this result holds true for all $k$-equal arrangements,
$k\geq 3$ (see \cite{perspectives}). That is
$$A^{n-k}_1(\Delta(B_n)) \simeq \pi_1(M_{n,k}).$$

So, if we want to compute the first Betti number, $H^1(M_{n,k})$
of the complement (in ${\mathbb R}^n $) of the $k$-equal arrangement, it will suffice to
compute the rank $A^{n-k}_1(\Delta(B_n))^{ab},$ which means counting
the number of equivalence classes of $6$-cycles in $\Gamma^{n-k}_{\Delta(B_n)}.$
The vertices of this graph correspond to the maximal chains of $\Delta(B_n)$
and there is an edge between two such maximal chains if they share (when viewed
as simplices of $\Delta(B_n)$) at least an $(n-k)$-face.
Thus $\Gamma^{n-3}_{\Delta(B_n)} (=\Gamma(B_n))$
is a subgraph of $\Gamma^{n-k}_{\Delta(B_n)}$
for all $k\geq 4,$ with both graphs having the same set of vertices.
In fact,
$$\Gamma^{n-3}_{\Delta(B_n)} \subset \Gamma^{n-4}_{\Delta(B_n)}
\subset \Gamma^{n-5}_{\Delta(B_n)} \subset \ \cdots ,$$

and the only difference between any two of these graphs is that
$\Gamma^{n-k}_{\Delta(B_n)}$ has more edges than $\Gamma^{n-j}_{\Delta(B_n)}$
for all $k>j.$ Already when $k=4$ it is easy to see that the graph
$\Gamma^{n-4}_{\Delta(B_n)}$ no longer has {\it primitive} 6-cycles.
Indeed, if $\sigma$ is a permutation of a $6$-cycle, and $m$ is the corresponding
maximal chain of $B_n$, then the maximal chain
$m'$ corresponding to the permutation $\sigma s_is_{i+1}s_i$ differs
from $m$ in exactly two levels. Hence, the permutations $\sigma$ and
$\sigma s_is_{i+1}s_i$ are adjacent in $\Gamma^{n-4}_{\Delta(B_n)}$
as their corresponding simplices share an $(n-4)$-face. Thus,
in $\Gamma^{n-4}_{\Delta(B_n)}$ all $6$-cycles are concatenation of at least two $4$-cycles,
which means that every cycle is contractible and $A^{n-4}_1(\Delta(B_n))$ is the trivial
group. Moreover, given the fact that $\Gamma^{n-k}_{\Delta(B_n)}$ contains
$\Gamma^{n-4}_{\Delta(B_n)}$ for all $k\geq 4$ none of these graphs
have primitive cycles of length $\geq 4,$ and we have the following theorem.

\begin{theorem}
For all $k\geq 4,$ $A^{n-k}_1(\Delta(B_n))=1.$
\end{theorem}

Thus, we can conclude that for all $k\geq 4,$ the first Betti number,
$H^1(M_{n,k}),$ of the complement of the $k$-equal arrangement is
zero, a fact that can also be recovered from Theorem $1.1 (b)$
of \cite{BW}.

\textbf{Conclusion}. As we mentioned in the introduction, once all is said and done
we believe that the crux of the paper lies in
the combinatorial techniques that were introduced.
In particular, realizing that the $\Gamma$-graph of the order complex of a product
of two posets is obtained by taking the box product of three graphs, one of them being the new
{\it shuffle} graph, and removing some (easily identifiable) edges,
should prove useful in computing the discrete fundamental group of other posets.
Moreover, the translation of the grid-homotopy equivalence relations in the language
of equivalence of words (over the symmetric group), is also promising.
Indeed, all Coxeter arrangements will have such algebraic structure that could prove
fundamental in computing the $A_1^q$-group of the appropriate complexes.

The authors would like to thanks the referees for their helpful comments
and advice.

\end{document}